\documentclass[citeauthoryear]{article}

\usepackage{amsmath}
\usepackage{amsfonts}
\usepackage{amssymb}
\usepackage{graphics}
\usepackage{latexsym}
\usepackage[pdftex]{graphicx}
\usepackage[mathscr]{eucal}
\usepackage{mathrsfs}
\usepackage{amsfonts}
\usepackage{amsthm}
\usepackage{amsmath}
\usepackage{dcolumn}
\usepackage{tabularx}
\usepackage{tabulary}
\usepackage{lscape}
\usepackage{multirow}
\usepackage{hhline}
\usepackage{amssymb}
\usepackage{amsfonts}
\usepackage{booktabs}
\usepackage{xcolor}
\usepackage{float}
\usepackage{ushort}
\usepackage{setspace}
\usepackage[utf8]{inputenc}
\usepackage{hyperref}
\usepackage{csquotes}
\usepackage[T1]{fontenc}
\usepackage[utf8]{inputenc}
\usepackage{textcomp}
\usepackage{eurosym}
\newtheorem{corollary}{Corollary}

\newtheorem{proposition}{Proposition}

\begin{document}

\title{\thispagestyle{empty} Optimal Lockdown Policies driven by Socioeconomic Costs\thanks{\textbf{Acknowledgments.} We benefited from discussions and comments from Pol Antras, Costas Azariadis, Giorgio Calcagnini, Samuel Bowles, Gian Italo Bischi, Laura Gardini, Germana Giombini, Giuseppe Travaglini, Davide Ticchi,   Harald Uhlig, Georges Zaccour, Charles Wyplosz. }}

\author{Elena Gubar\thanks{\textit{Corresponding author.} Faculty of Applied Mathematics   and Control Processes, St. Petersburg State University, 7/9 Universitetskaya nab., St. Petersburg, 199034 Russia. E-mail: e.gubar@spbu.ru}  \and Laura Policardo\thanks{The Customs and Monopolies Agency (Agenzia delle Dogane e dei Monopoli), Sesto Fiorentino, Italy. E-mail: laura.policardo@gmail.com  } \and Edgar J. S\'{a}nchez Carrera\thanks{Department of Economics, Society and Politics at the University of Urbino Carlo Bo, Italy. Research Fellow at CIMA from the Autonomous University of Coahuila, Mexico. E-mail: edgar.sanchezcarrera@uniurb.it} \and Vladislav Taynitskiy \thanks{Faculty of Applied Mathematics   and Control Processes, St. Petersburg State University, 7/9 Universitetskaya nab., St. Petersburg, 199034 Russia. E-mail: tainitsky@gmail.com}}

\maketitle

\begin{abstract}
\noindent
In this research paper we modify a classical SIR model to better adapt to the dynamics of COVID-19, that is we propose the heterogeneous SQAIRD model where COVID-19 spreads over a population of economic agents, namely: the elderly, adults and young people.  We then compute and simulate an optimal control problem faced by a Government, where its objective is to minimize the costs generated by the pandemics using as control a compulsory quarantine measure (that is, a lockdown). We first analyze the problem from a theoretical perspective, claiming that different lockdown policies (total lockdown, no lockdown or partial lockdown) may justified by different cost (concave or convex) structures of the economies. We then focus on a particular cost structure (convex costs) and we simulate a targeted optimal policy vs. a uniform optimal policy, by dividing the whole population in three demographic groups (young, adults and old). We also simulate the dynamic of the pandemic with no policy implemented. Simulations highlighted the fact that: a) a policy of lockdown is always better than the \emph{laissez faire} policy, because it limits the costs that the pandemic generates in an uncontrolled situation; b) a targeted policy based on age of the individuals outperforms a uniform policy in terms of costs that it generates, being a targeted policy less costly and equally effective in the control of the pandemic.

\bigskip


\end{abstract}

\newpage

\section{Introduction}
This research paper develops a SQAIRD (Susceptible-Quarantined-Asymptomatic-Infected-Recovered-Dead) model for three groups of people: elderly, adults and young people. It studies an optimal lockdown problem where the planner wants to minimize a social cost function allowing for either convex or concave costs. Two sets of results emerge: (i) a policy of lockdown is always better than the laissez faire policy and (ii) targeted policies are better than uniform policies. Particularly we show that: (i) Endogenous lockdown with a starting and end time is imposed only when the cost of shifting individuals into quarantine is convex in the rate of quarantining. For a concave cost, there is a boundary solution with either no lockdown or complete lockdown. (ii) It typically pays to impose selective lockdown, focusing on the group that is ``least costly'' in being locked down. The results are derived theoretically and illustrated by way of numerical simulation. While there are already several papers allowing for pandemic dynamics that are richer than an SIR model. Two examples with population of different age groups and also with dynamic modeling are Brotherhood et al (2020a,b) and Glover et al (2020). Although they don't study exactly a SQAIRD model, they study the implications where lockdown policy is always better than no-policy and targeted policies are also better than non-targeted. Our approach is quite similar to the one developed by Acemoglu et al. (2020), from which they conclude that the isolation of individuals belonging to particular social groups is an effective means to limit the number of deaths and the incidence of COVID-19 and the economic loss of the virus. These authors state, based on their model, that a policy of isolation of the elderly from other age groups (which they describe as ``group distancing'') is a very efficient means of addressing the trade-off boundary between incidence of COVID-19 and the loss of production due to the disease. Instead, in our approach we show that it is not the isolation of an age group but the separation into age groups that results in an optimal quarantine policy with the objective of minimizing different structures of socioeconomic costs. Similarly, the modeling of Acemoglu et al., only focuses on the economic impact of death rates and does not specifically take into account the economic losses that occur when people are affected by the virus but do not die. And this fact is addressed in the present research article.

The recent pandemic of COVID-19 that spread all over the world has put Governments at hard testing, because they have to manage a global health crisis with dramatic effects on both human lifes and economies. At the time this paper was written, indeed, ILO estimated a loss of nearly 25 million jobs worldwide due to the pandemic, and in the OECD area, economic losses of about 10\% for several European countries (among those, Italy, France and United Kingdom) which, moreover, are estimated to grow in the near future (Petrosky-Nadeau and Valletta, 2020).
However, policies undertaken by Governments aimed at managing this crisis have been very different. Some Governments preferred not to take any direct action to limit the spread of the virus (the case of Sweden, for example, which decided not to undertake any lockdown in order to limit the losses coming from their economic activities)\footnote{See Claeson and Hanson, 2020: https://www.ncbi.nlm.nih.gov/pmc/articles/PMC7755568/}, while others have preferred a full lockdown even before recording a single death as a cause of the virus, that is the case of Finland (Moisio, 2020).

One may wonder if there are any reasons for such differences in managing the same problem, differences that apparently are not justified by different culture and lifestyles, being undertaken by rich and developed countries, moreover endowed with similar cultural backgrounds and geological/climatic characteristics. In this paper we explore, from an optimal control point of view, what might be the reasons behind such different policies undertaken by Governments, and moreover, we estimate the optimal duration and the estimated costs of a lockdown policy under the hypothesis of convex social cost functions. We explore this optimal policy by comparing an uniform lockdown policy with a targeted one on the basis of demographic differences (i.e. age groups). If a policy of lockdown or no lockdown can be justified by convex/concave social cost functions, targeted lockdown policies may be undertaken in case of convex social cost functions in order to minimize the economic costs of the pandemic.

To the best of our knowledge, a targeted lockdown policy has not been undertaken by the majority of countries who suffered the consequences of the pandemics.
We believe that exploiting the different impact that the virus has on different demographic groups, a country may reduce the economic losses of the pandemics by imposing differentiated lockdowns for people belonging to different demographic groups.

Literature on the epidemiology of COVID-19 claims clearly that the risk of serious health complications vary greatly among different demographic groups. Mortality rates of Covid-19 for people aged 65+ are indeed 60 times higher than those of people aged 20-49 (Yanez, N.D., Weiss, N.S., Romand, JA. et al., 2020), and just this fact can be worth examining the possible benefits of targeted lockdown policies based on age in terms of social costs. Intuitevely, indeed, letting some (possibly the weaker) individuals to avoid contacts with others through a stricter and longer lockdown policy allows lower risk groups (the younger) to obey a shorter and less strict lockdown, thus reducing the indirect costs of management of the infection (i.e. the costs associated to losses of GDP due to the interruption of economic activities).

Our baseline model represents a population divided in five sub-populations: Susceptible (S), Quarantined (Q), Asymptomatic (A), Infected (I), Recovered (R), and Dead (D), (in the following, ``SQAIRD'' model). Each sub-group is moreover stratified into three groups based on age (youngs, adults and old people) and each group has different risks of death, health complications and so on. We set up an optimal control problem where the Government wants to minimize a social cost function subject to the equations determining the dynamics of the subgroups of population as a consequence of the spread of the virus.
As a controls, the government has to decide how many people belonging the the susceptible group have to be isolated, for each of the subgroups of young, adults and old.

\bigskip

The reminder of the paper is organized as follows: section 2 introduces a literature review on the topic; section 3 illustrates the epidemic process; section 4 presents the theoretical optimal control problem, and derives the conditions under which is optimal a total lockdown, no lockdown at all or a mixed policy; section 5 calibrates the model under the hypothesis of convex cost functions, defines the functional form of the objective and state equations and illustrates the numerical results. It is important to point out that calibration of the costs and parameters used in the model is based on data available for Italy. Section 6 concludes.

\section{Literature Review}

Since the spread of the COVID-19 pandemic, scholars have devoted an unstoppable and pressing effort in studying both the economic impact of the pandemic and the policy effects aimed at stemming the problem of infection on the economies.

The economic literature that is available on the topic at the time this paper is written can be classified roughly into two branches; the first one aimed at studying the impact of the pandemic on several social and economic indicators, mainly from an empirical point of view, and the second one which makes use of mathematical models and, in second instance, optimal control theory to compare and forecast the effect of a policy under different scenarios. Both methods have advantages and pitfalls. While the first method indeed allows to compare completely different policies that act on different factors affecting the spread of the virus on a deterministic time span, the second one allows to compare the same policies at most on different subgroups of population, as compared to the same policy applied to the whole population. In order to compare different policies, indeed, in an optimal control framework, several problems have to be solved. Another advantage of the first method is that it does not require assumptions(that might be wrong) on the parameters of the model, because it is often based on real data, while the second approach requires assumption of the dynamics of the spread of the virus and other parameters that might not respect completely the reality of the facts.  Among the pitfalls of the first methods, indeed, we account for the possible biases due to the omitted variable problem and simultaneity that pervade almost all the econometric studies in all disciplines. The second methods, indeed, does not suffer of this problem because the effect of the policy under study is already purified of this potential distorting effect, and moreover, it allows to compare ``all things being set equal'' targeted policies over different subgroups of population.

\bigskip
In the first group of literature we can number many recent papers. The majority of the papers concerning the outbreak of COVID-19 pandemic have tried to evaluate different policies from an empirical perspective. With regard to the effectiveness of the policies undertaken by governments aimed at stemming the infection rate, Zhixiam Lin and Meissner (2020)  found that the policies most strongly and statistically significantly associated with slowing the growth rate cases were public transport closures, enforced workplace closures, limited domestic travel, and restrictions on international travel. School closures and limits on public events were not found to be statistically significant. For example, from a cross-country perspective, the effectiveness of non-parmaceutical policies are found to be negatively correlated with percapita GDP, population density and surface area, while they are positively correlated with health expenditure and fraction of physicians in the population (Bargain and Aminjonov, 2020; Castex et al., 2020). Other scholars investigated the difference between different policies that a government could implement to limit the losses of the pandemic (both human and economic), from other points of view: Brotherhood et al. (2020a,b) investigated, from a cost-benefit analysis, the implementation of two policies: shelter-at-home order, requiring that individual stay at home at least 90\% of their time, and test-and-quarantine policy, which imply an additional cost due to testing, but only people who are found infected are isolated. The first policy is assumed to last 26 weeks and it cuts dramatically the labor supply of the young, with a substantial decrease in GDP. The second policy, conversely, requires infected agents to stay longer isolated, even if this is against their best interest, but has the advantage of decreasing the number of unaware infected around with a smaller decrease of GDP. In this vein, we have assumed a deterministic duration of the policy intervention, while they opened up to the possibility to implement targeted policies over different demographic groups. Other scholars indeed offer different solutions to the spread of the virus: Fajgelbaum et al. (2020) analyze targeted lockdown policies from a spatial point of view. They claim indeed that locking down more densely populated areas, while leaving others in business, may help reducing the spread of the virus minimizing the losses in GDP.
This kind of works are numerous and one of the main problem is that the estimations obtained might be biased due to the policies undertaken by government that might have mitigated the effects.

\bigskip

In the second branch of literature, we can number several attempts to study the effect of a policy under different scenarios. From a comparative static approach, we encounter Bandyopadhyay et al. (2021) who characterize the timing, intensity and duration of a lockdown by balancing the opposite effects that a lockdown generates, that is to say, the beneficial habit formation like social distancing and personal hygiene, and the economic costs, which are substantial. From a dynamical perspective instead, we number several attempts to perform an optimal policy analysis.
Scholars indeed have started incorporating economic trade-offs and conducting optimal control policy analysis within the SIR framework, undertaking an optimal control analysis in single and multiple-group models (starting with the early related contributions by Geoffard and Philipson, 1996; Fenichel et al., 2011; Fenichel, 2013; Goyal and Vigier, 2015; and recently: Goenka and Lio, 2019; Acemoglu et al., 2020; Rowthorn and Toxvaerd, 2020; Eichenbaum, Rebelo and Trabandt 2020; Alvarez, Argente and Lippi 2020; Jones, Philippon and Venkateswaran, 2020; Farboodi, Jarosch and Shimer, 2020; Garriga et al., 2020). All of those research, however, do not consider both: an extension of the standar SIR modeling (which is conisdered here, as explained in what follows, i.e. a SQAIRD dynamic model), and the role and type of the socioeconomic costs.

Rowson et al. (2020), for example, consider an optimal lockdown policy with progressive release of population under isolation once that the prevalence of infected people decreases beyond a given threshold (determined by the capacity of the health system). These authors however do not study targeted policies based on demographic characteristics of people, and therefore do not exploit fully the benefits given by the fact that the younger have a greater resistance against the infection, so the release of people could generate more gains in economic terms if this is done with a specific target of population.

Piguillem and Shi (2020) compare - from an optimal control perspective - different policy interventions aimed at reducing the economic consequences of the pandemic. As in many recent contributions, they use a SIR model to identify the optimal Government response to the spread of the virus, which depends on the curvature of the welfare function considered and maximized, and on the value of a life.
Their main result is that a suppression strategy (i.e. total lockdown) is never optimal if the virus cannot be eliminated from the population, but, depending on the value of a life (which threshold is 13 years of annual income, according to the authors), mitigation, that is to say, mild interventions aimed at ``flattering the curve of infection'' or no intervention can be optimal.
If testing is feasible, it is possible to reduce the duration of the policies, because the asymptomatic infected people can be isolated and therefore cannot be able to infect other susceptible.
As in other papers, the authors fail to exploit the demographic differences of the population to reduce further the losses caused by the lockdown, total or partial.

An effort to measure the benefits of targeted lockdown policies with respect to uniform ones has been done by Acemoglu et al. (2020). Using an heterogenous SIR model, they find that optimal policies targeting different risk/age groups outperform uniform policies, and most of the gains can be realized by imposing stricter lockdown policies for the most vulnerable group, thus allowing the other groups to obey a less strict lockdown.
Our paper is, somehow, similar in several aspects to that of Acemoglu, with the exception that we extend the analysis to the cases where the optimal policy implies full lockdown or no lockdown at all. In addition, we formally demonstrate the existence of a solution to the optimal control of the coronavirus over a multilayer network population of economic agents, and that they are affected by different socioeconomic costs, whether concave or convex (for example, economies of scale).
Therefore our research paper has some important  advantages, i.e.: i)  we develop the SQAIRD model and its optimal control problem with different groups of economic agents and targeted policies, ii) we show the solutions  this problem with different cost structures (convex or concave), and iii) numerical simulations corroborate the optimal lengh of a targeted lockdown policy in an environment of different economic agents driven by socioeconomic costs, and the possibility to analyze the asymptomatic infected (this is assumed in the choice of the parameters in the simulation).


\section{The epidemic process}\label{epidemic}

We model the spread of the virus in a population by extending a classical Susceptible-Infected-Recovered-Dead (SIRD) model by Allan (2008), Altman et al. (2011) Gubar and Zhu (2013), and Gubar et al. (2017). Within this model, we first divide the whole population into three demographic groups, denoted by $p=1,2,3$, which represent, respectively, youngs (1), adults (2), and the elderly (3), and  Within each demographic group, the population is divided in further five subgroups, let's say the \textit{Susceptible} $(S)$, \textit{Quarantined} $(Q)$,  \textit{Asymptomatic} $(A)$, \textit{Infected }$(I)$, \textit{Recovered} $(R)$, and \textit{Dead} $(D)$. This division of the population, as previously mentioned, extends the classical SIRD model by introducing two additional groups, that is, the Quarantined and Asymptomatic, therefore letting become the SIRD model a ``SQAIRD'' model.

The group of Susceptibles are those who are exposed to the virus (because they do not have immunity) and are, at the beginning of the pandemic, the majority of the population. The group indicated by $Q$ is the group quarantined (i.e. isolated at home), by law or by their will. The group $A$ is instead a group of people who catched the virus, but did not show evident symptoms. These people can behave like a susceptible but they can infect other people as well. Once infected, a person shift in the group of Asmptomatic, then he can either recover (then shift into the group of Recovered $R$) or develop symptoms, then shifting in the group $I$ of the Infected with symptoms with probability $\alpha_p \in (0,1)$, for each demographic group $p=1,2,3$. An Infected, instead, can either recover and shift into $R$ or die (then shift into the $D$ group).
Once a person shift into the $R$ or $D$ groups, it cannot be infected anymore (we assume that a person that recovers, develops immunity and cannot be reinfected).

At each time $t$ and for each demographic group $p$, a fraction of Susceptible, Quarantined, Infected, Recovered and Dead exist, and are denoted by $S_p(t)$, $Q_p(t)$, $A_p(t)$, $I_p(t)$, $R_p(t)$ and $D_p(t)$, respectively.

By normalization:
\begin{equation}
   \sum_{p=1}^3 S_p(t)+ \sum_{p=1}^3 Q_p(t) +\sum_{p=1}^3 A_p(t) +\sum_{p=1}^3 I_p(t)+\sum_{p=1}^3 R_p(t)+\sum_{p=1}^3 D_p(t)=1,
\end{equation}
that is to say, the sum each fraction of Susceptible, Quarantined, Asymptomatic, Infected, Recovered, and Dead people for each group of individuals is equal to one. The whole population is assumed constant and equal to $N = 1$ during all the period of analysis.

At the beginning of the epidemic ($t=0$), the majority of the people belongs to the Susceptible group, and a small fraction of people are infected  (see Sahneh et al. (2013), Mieghem (2009), Taynitskiy et al. (2017)). Hence, for each demographic group $p$, the initial states are:
\begin{eqnarray}
    \nonumber 0< &S_p(0)& =S_p^0 \\
    \nonumber 0 \leq &Q_p(0) &=Q_p^0 \\
    \nonumber 0 \leq &A_p(0) &=A_p^0 \\
    \nonumber 0< &I_p(0) & =I_p^0 \\
    \nonumber 0 \leq &R_p(0) &=R_p^0 \\
    \nonumber  N_p-S_p^0- Q_p^0-I_p^0-R_p^0 =&  D_p(0) & =D_p^0
\end{eqnarray}
where $\sum_{p=1}^3 N_p = 1$ and $S_p(t)+Q_p(t)+A_p(t)+I_p(t)+R_p(t)+D_p(t)=N_p$. Figure 1 illustrates the epidemic process in this SQAIRD model. In such a scheme, the arrows indicate the infection process among groups, ending at the node of recovered, quarantine or dead agents.

\begin{figure}[H]
\begin{center}
\includegraphics[width=0.5\linewidth]{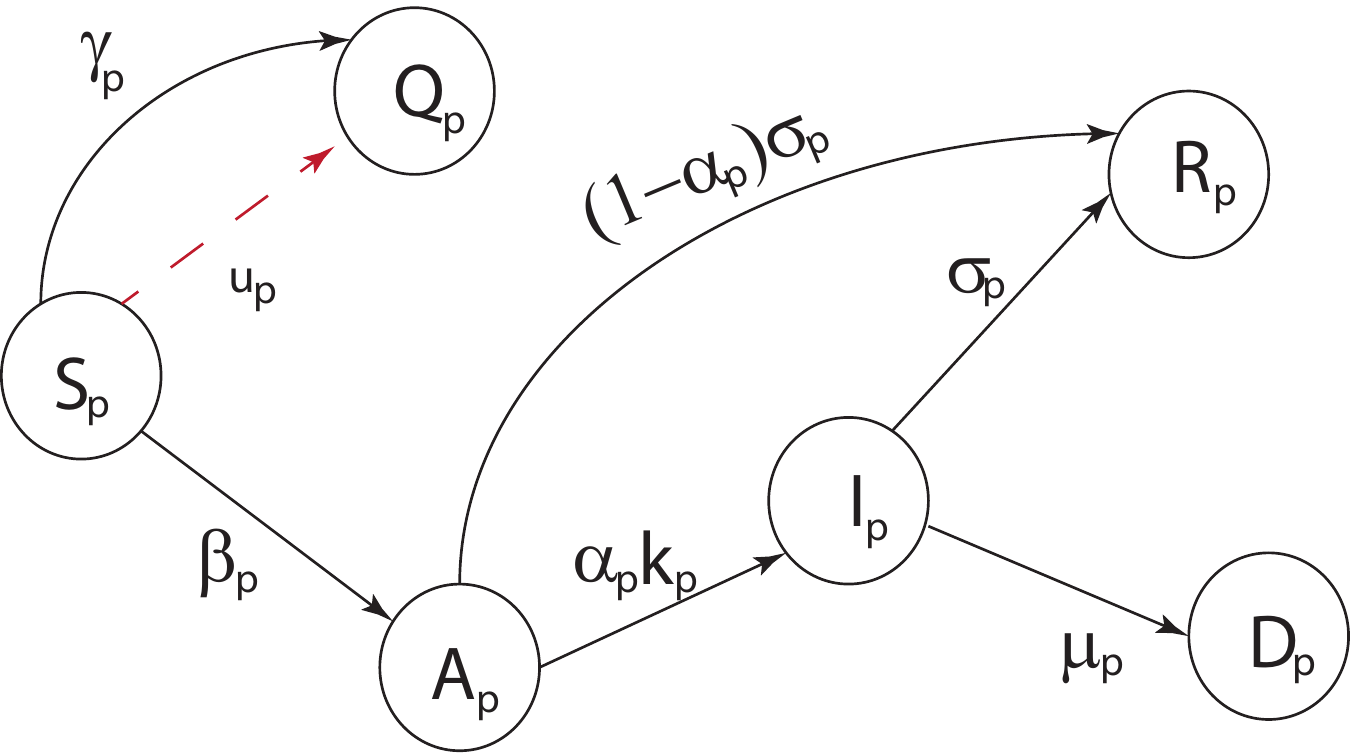}
\label{SQEIRD}
\caption{The scheme of epidemic process, index is $p=1,2,3$. Nodes corresponds to the fraction of infected people in the entire population.}
\end{center}
\end{figure}

Figure 1 clearly defines the process of epidemic: a Susceptible individual would be infected and shift into the category of Asymptomatic with probability $\beta$, which we assume constant through all the three categories of individuals (young, adults and old people). At this point, he can recover naturally with probability $(1-\alpha_p)\sigma_p$ ($\alpha_p$ represents the probability to develop symptoms once infected and $\sigma_p$ is the natural rate of recovery) or develop symptoms with probability $\alpha_p$. The parameter $k_p$ represents the ``speed'' at which an Asymptomatic shifts into the Infected group. Once infected, he can recover with probability $\sigma_p$ or die with probability $\mu_p$. Notice that the probability to recover or to die differs between demographic groups, contrary to the probability of catching the virus (i.e. transiting into the group of asymptomatic).

We make the following simplifying assumptions:
\begin{itemize}
    \item[-] All the asymptomatic that will develop symptoms stay on average without symptoms for 5 days before shifting into the Infected group\footnote{According to the Report of the WHO-China Joint Mission on Coronavirus Disease 2019 (COVID-19), available at \url{https://www.who.int/docs/default-source/coronaviruse/who-china-joint-mission-on-covid-19-final-report.pdf}, symptoms are developed, on average, 5 days after infection.}. This value is embodied in the parameter $k_p$.
    \item[-] Asymptomatic people recover naturally with a rate $\sigma_p$. We make moreover the assumption that all the Infected people (that is to say, those who develop symptoms) are treated with medical procedures, which vary accorging to the severity of the infection. Once infected and treated, an individual may recover with probability $\sigma_p$ or die with probability $\mu_p$.
    \item[-] Both quarantine measures for the susceptibles and medical treatment (irrespective of whether this takes the form of hospitalization or isolation) for the infected have the effect of reducing the number of Susceptibles that are exposed to the infection and the number of Infected that Susceptibles can get in touch with, so have the effect of slowing the spread of the virus. Asymptomatic persons can still be able to infect other susceptibles.
    \item[-] A Susceptible person can also shift into the Quarantined group by his own free will, she/he decides with probability $\gamma_p$ to isolate himself (i.e. stay at home). Another reason why a susceptible person can shift into the quarantined group is by law (the Government undertake a lockdown policy) with probability $u_{p}$. These measures (irrespective of whether they are defined by people's will or by law) have the effect of reducing the number of susceptible people that can be infected, so these measures are also able to reduce the speed of circulation of the virus.
\end{itemize}

According to the process of spread of the virus just described, we formalize it through the following differential equations:

\begin{eqnarray}
\frac{dS_p(t)}{dt} &=& -\beta S_p(t)(A(t)+I(t))-(u_{p}(t)+\gamma_p)S_p(t),  \label{deltaS} \\
\frac{dQ_p(t)}{dt} &=&  (u_{p}(t)+\gamma_p)S_p(t), \label{deltaQ} \\
\frac{dA_p(t)}{dt} &=&  \beta S_p(t)(A(t)+I(t)) - \alpha_pk_pA_p(t)-(1-\alpha_p)\sigma_pA_p(t),\label{deltaA} \\
\frac{dI_p(t)}{dt} &=& \alpha_pk_pA_p(t)-( \sigma_p+\mu_p)I_p(t), \label{deltaI} \\
\frac{dR_p(t)}{dt} &=&  \sigma_p I_p(t)+(1-\alpha_p)\sigma_p A_p(t), \label{deltaR} \\
\frac{dD_p(t)}{dt} &=& \mu_p I_p(t). \label{deltaD}
\end{eqnarray}
In equations \ref{deltaS} - \ref{deltaD} the terms $A(t)=\sum_{p=1}^{3} A_p(t)$ and $I(t)=\sum_{p=1}^{3} I_p(t)$ indicate the numerosity of the Asymptomatic and Infected in the society at time $t$.

\bigskip

Equation \ref{deltaS} describes the behaviour through time of the Susceptible group ($S_p$). During the pandemic this group is intended to decrease as time runs. A fraction of Susceptibles indeed gets the infection with probability $\beta$ each time a contact with an infected or asymptomatic person happens. The probability of getting infected is indeed proportional to the stock of infected people freely circulating (basically the virus spread as a cause of infected individual who get the infection from 1 to 5 days before) and it is equal to $\beta (I(t)+A(t))$, while another reason why this group decreases is that individuals shift in the quarantined group because either they get isolated themselves by their will or because the government force them to stay home. This probability is given by $(u_{p}(t)+\gamma_p)$. So, the fraction of Susceptibles flowing from this group to others (in particular, to the Asymptomatic or Quarantined) is described by Equation \ref{deltaS}.

The group of quarantined - as we already mentioned - represents the fraction of susceptible under compulsory isolation ($u_{p}$) or self-isolation ($\gamma_p$). The dynamics denoting the behaviour of group $Q_p$ is denoted by Equation \ref{deltaQ}.

Equation \ref{deltaA} describes the dynamics of group of Asymptomatic ($A_p$). With the probability $\alpha_p$, an Asymptomatic person can develop symptoms and become Infected. This happens with the rate $k_p$, reflecting the fact that an asymptomatic person usually takes 5 days to develop symptoms. With the probability $(1-\alpha_p)\sigma_p$ the asymptomatic person will recover from the virus without developing any symptom.

Equation \ref{deltaI} describes how the Infected group ($I_p$) behaves through time. On one hand, indeed, there is a positive contribution to this group from the asyntomatics, which become infected (i.e., develop symptoms) at a rate $\alpha_p k_p$. On the other hand, this group decreases because some of the infected recover and some die. This counterbalancing forces will lead the curve describing the amount of the infected during the pandemic be an inverse U, with an increasing trend in the first periods of the pandemic, and a negative trend thereafter.

As previously mentioned, an infected individual would recover with probability $\sigma_p$ and we assume that an infected person gets detected after 5 days of the first appearance of symptoms.\footnote{\url{https://www.cdc.gov/coronavirus/2019-ncov/symptoms-testing/symptoms.html},\\ \url{https://www.healthdirect.gov.au/coronavirus-covid-19-symptom-faqs}} Recovery happens, on average, after 30 to 40 days from infection (we assume 30 days for young, 40 days for adults and 60 days for the elderly).  Infected individuals that do not recover, are assumed to die with probability $\mu_p$. So, while the group of Infected increases because the Susceptibles get sick with probability $\beta (I(t)+A(t))$, this group decreases because the infected can either recover or die with probability $(\sigma_p +\mu_p)$.

Equation \ref{deltaR} represents the dynamics of the Recovered group ($R_p$), which increases through time  because a fraction of Infected defeats the virus (this may be due to effective health treatments or maybe naturally), and this happens with probability $\sigma_p$. This probability determines the positive flow of individual coming from the Infected and Asymptomatic groups. Equation \ref{deltaR} formalises this dynamics.

Finally, Equation \ref{deltaD} represents the dynamics of the group of deaths ($D_p$). During the pandemic, it increases - obviously - through time according to the fraction of Infected who cannot recover. Then, in each period of time, this group increases by a proportion of $I_p(t)$ equal to $\mu_p$.

In the following sections, we introduce the possibility that the government undertakes a policy through compulsory lockdown of the Susceptibles, which we model by $u_p$, for each demographic group. Depending on the various costs structures of the economy, different scenarios may appear that can justify different policies undertaken by similar countries (both in terms of level of infection and other characteristics). Moreover, it is possible to show, with numerical simulations, that targeted policies over the different demographic groups allow to reduce the economic losses of the lockdown.

\newpage

\section{The optimal control problem: a theoretical perspective.}\label{costbenefits}
The spread of a potentially lethal virus may have devastating effects in an economy. First of all, because ill and lost lives produce a huge impact on the labor market due to the temporary and permanent loss of human capital. A negative shock in the labor force negatively affects production and therefore may impoverish families and decreases government's entrance. The pandemic has also a strong negative effect on the health system, which has to bear the treatment costs of a great number of infected, often without the necessary technical and human resources.
On the political side of the problem, a pandemic that is not well managed may affects the political orientation of the citizens, so the future re-election of the political party in charge may be at risk. So, a correct management of the emergency caused by the pandemic is of vital importance for any Government. Therefore, there is the need to find an effective way to control the COVID-19 pandemic with minimal economic and social disruptions.

Optimal control theory, in this special case, may help understand  how to apply one or more time-varying control policies to a nonlinear dynamical system, in such a way that a given objective function (that in this case could be a cost and profit function) is optimized (Kissler et al., 2020; Perkins and Espana, 2020). Here, we apply optimal control theory to determine optimal strategies for the implementation of lockdown policies to control a pandemic. As we will see shortly, an optimal strategy (i.e. a strategy that minimizes the costs of a lockdown) may depend on the socioeconomic cost structure of the economy, in particular, on the concavity/convexity of these functions.

Let us introduce the objective function (i.e. the function that has to be minimized or maximized, depending on the context) that the government has to optimize using as controls $u_p(t)$, $p=1,2,3$.
This objective function is determined by the sum of direct and indirect costs related to the management of the pandemic, eventually summed to a ``profit function'' concerning those who re-enter the job market after recovery.
In other words, the Government has to optimize the balance between profits and costs through the application of a policy, which consists in the obligation for a fraction of the Susceptibles to stay isolated and do not keep in touch with other individuals. This policy is differentiated between demographic groups, that is to say that the obligation for the young susceptibles could be different than the obligation for the adults and old susceptibles.
Let us define the following:
\begin{itemize}
    \item[-] $f_p(I_p(t))$ is a function measuring the direct costs of infection, with a non-decreasing, twice differentiable function, such that $f_p(0)=0$, $f_p(I_p(t))>0$ for $I_p(t)>0$. These costs can be interpreted by the price the Government pay to allow the health system treat the infected and increase the probability of recovery for the treated individuals. The infected that are not detected because they have not symptoms are not treated and therefore there are not direct costs for them.

    \item[-] $g_p(R_p)$ represents the profit function for people that are recovered from the illness. This function is assumed to be non-decreasing and differentiable.


    \item[-]Imposing quarantine for Susceptibles or in the case that some individuals have decided to self isolate, so they cannot have contacts with infected people, is also costly: while there is no direct cost of medical treatment, because they are not sick, there is the fact that those people stay away from the job market, meaning that they stop producing for the whole period of isolation. Similarly for the infected, a production that is interrupted for a period of time implies an impoverishment of the people (and this may call for a corrective action in support of the necessary consumption for the people heavily affected by the lockdown) while the revenues of the government decreases. We model such socio-economic costs by using functions $h_1(u_1(t))$, $h_2(u_2(t))$ and $h_3(u_3(t))$, that depend on the fraction of isolated Susceptible belonging to the young (1), adults (2), and old (3) groups. These functions $h_{p}(u_{p}(t))$ are assumed to be increasing in the arguments $u_{p}(t)$, and twice differentiable such that $h_{p}(0)=0$, $h_{p}(u_{p}(t))>0$ when $u_{p}(t)>0$ for each $p=1,2,3$.

    \item[-] Finally, each life lost produces a cost for the society. We do not model it in the theoretical study of the problem because, as it will be clear later, this cost is assumed to be constant for each individual in a given demographic group (but it varies between demographic groups).

\end{itemize}

Hence, the problem that an hypothetical Government must face becomes then to maximise the difference between benefits and costs, using as controls the fraction of susceptible isolated $u_p(t)$, for each demographic group $p=1,2,3$ and taking into account the dynamics of the pandemic which is represented by equations \ref{deltaS})-(\ref{deltaD}. The cost-profit function is mathematically represented by the following equation:
\begin{equation}\label{functional_J}
J(u_1,u_2,u_3)=J_1(u_1)+J_2(u_2)+J_3(u_3)
\end{equation}

\noindent where

\begin{equation}\label{functional_Jp}
J_p(u_p)=\displaystyle\int_{0}^{T} \Big[  f_p(I_p(t))+h_{p}(u_{p}(t))-g_p(R_p)\Big]dt, \quad p=1,2,3.
\end{equation}

\bigskip

Therefore, the problem that the Government faces is to optimize the objective function represented by Equation (\ref{functional_Jp}), subject to the SQAIRD dynamics represented by equations (\ref{deltaS})-(\ref{deltaD}).

Recall that we define the costs $h_p(u_p(t))$, $p=1,2,3$ as twice differentiable and increasing in his argument, but we did not specify the sign of the second derivative, which may imply convexity or concavity. Such behaviour of the functions indeed changes dramatically the type of optimal duration that a policy of lockdown may take.


\begin{proposition}
If the cost functions $h_p(u_p(t))$, $i=1,2,3$ are concave, that is, the second derivative with respect to its argument is negative or equal to zero, then there exist an optimal $t_0 \in [0,T]$ such that for any $p=1,2,3$,
$$u^*_p(t)=\left\{\begin{array}{l} u_{max},\ \mbox{for}\ 0\le t\le t_0;\\
 0,\hskip 20pt \mbox{for}\ t_0<t\le T.\end{array}\right.$$

\bigskip

\noindent This means that the solution is a corner solution, implying that there exists a threshold period $t_0$ such that below such a threshold, the optimal policy is to impose total lockdown ($u_{max}=1-\gamma_p$), for the whole population, and above this threshold, the optimal policy is not to impose any lockdown. The reason for such result is quite intuitive. As long as the period of lockdown is short enough, the economic losses in terms of production are smaller relative to the benefits due to the missed infections, whose treatment is costly. If the period of lockdown increases above the threshold, then it is preferable to bear the costs of infection rather than the losses due to the missed production.

\bigskip

\noindent If, rather, the costs $h_p(u_p(t))$, $p=1,2,3$ are convex, that is to say, the second derivative with respect to its argument is greater or equal to zero, there there exists two time moments, which we denote by $(t_0, t_1)\in [0,T]$ such that for any $p=1,2,3$ and $\zeta(t) \in(0,u_{max})$,
$$
u^*_p(t)= \left\{\begin{array}{l}u_{max}, \hskip 6pt 0\le t\le t_0;\\
 \zeta(t),\hskip 12pt t_0<t\le t_1;\\
 0,\hskip 23pt t_1<t\le T.\\
\end{array}\right.$$

\bigskip

\noindent Note that if the costs of quarantine are convex, indeed, in addition to the corner solutions identified previously, there exists an interior solution identified by the fraction $\zeta(t)$ of isolated population if the period of lockdown lies between these two thresholds $t_0$ and $t_1$.

\end{proposition}

This model of optimal lockdown can indeed explain the dramatically different policies undertaken by several countries in Europe and in the World.
For example, it is not surprising indeed that countries like Sweden had decided not to enforce the citizens to any lockdown, limiting to suggest to the people the well know recommendations such as social distancing and personal hygiene. The opposite case of Sweden could be Finland, which imposed a strict lockdown for all the productive activities (with very few exceptions) and strict curfew for citizens even before recording a single death due to the virus. This, of course, has implied a total revolution in the way people work and behave, with a massive implementation of home working, and further investments in physical capital to substitute labor whenever possible.
These cases can be classified among those facing concave costs, typical of advanced/tertiary economies, where higher levels of physical and human capital make innovations and flexibility easier to digest for the economic system.

Economies based on manufacturing instead, or where the economy is driven by sectors - like tourism - which can be heavily damaged by lockdown policies, or where the level of government debt and unemployment are high, may face convex cost functions $h_p(u_p(t))$, which can be responsible for the initial hesitancy of Governments in implementing hard virus-restraint policies, like for example Italy, Spain and UK. Governments of these countries indeed waited in applying global confinements measures, limiting - in the first period - to confine only the areas where the incidence of the virus was higher, and progressively extended the lockdown to the whole nation.

\bigskip

We can state as corollaries that:
\begin{corollary}
An obvious consequence of the analysis conducted in this subsection is that at time $t_0\in (0, T)$, an optimal strategy for a government is to introduce a selective lockdown policy, which is aimed at protecting the population from the spread of the virus by selectively choose the individuals who are more vulnerable and/or less productive for lockdown, while allowing the others to keep their activity ``alive''.
\end{corollary}

\begin{corollary}
The corollary of the analysis indicates that as of moment $t_0\in (0, T)$, for the fight against COVID-19 it is really necessary to contemplate a strategy of greater flexibility for economic activities (but only if they strictly respect each of the rules for containing the virus) combined with more selective containment/lockdown measures of economic agents, aimed at protecting the most vulnerable groups.
\end{corollary}

In the following, we provide the formal proof of concave and convex cost function for the optimization of the Hamiltonian obtained from eq. (\ref{functional_Jp}) and Proposition 1. The reasoning offered is pretty intuitive, and allows the reader to figure out how the optimization performs in the two different cases.

\subsection{ Functions $h_p(\cdot)$ and the Hamiltonian are concave}
Let us define the Hamiltonian and the adjoint system for the initial system (\ref{deltaS}) - (\ref{deltaD}) of differential equations that describes the propagation of the virus (see Pontryagin (1962), Altman (2011), Gubar (2013)).  By  using  Pontryagin's maximum principle \cite{Pontryagin}, we  construct the optimal control $u(t)=(u_1(t), u_2(t), u_3(t))$ to the problem described above in Section \ref{epidemic}. To simplify the presentation, we use short-hand notations $S,I_1,u_1,$ etc. in place of $S(t),I_1(t),u_1(t),$ etc. Define the associated Hamiltonian $H$ and adjoint functions $\lambda_{S_p}(t)$, $\lambda_{Q_p}(t)$, $\lambda_{A_p}(t)$, $\lambda_{I_p}(t)$, $\lambda_{R_p}(t)$, and $\lambda_{D_p}(t)$, $p=1,2,3$ as follows:

\begin{equation}\label{Hamiltonian}
\begin{array}{l}
H=H_1+H_2+H_3,
\end{array}
\end{equation}
where
\begin{equation}\label{Hamiltonian_p}
\begin{array}{l}
H_p=f_p(I_p)+h_p(u_p)-g_p(R_p)+(\lambda_{A_p}-\lambda_{S_p})\beta S_p(\sum_{p=1}^{3} A_p+\sum_{p=1}^{3} I_p)+\\
\hskip 25pt (\lambda_{Q_p}-\lambda_{S_p})(u_p+\gamma_p)S_p+(\lambda_{I_p}-\lambda_{A_p})\alpha_p k_p A_p+(\lambda_{R_p}-\lambda_{A_p})(1-\alpha_p)\sigma_pA_p+\\
\hskip 25pt (\lambda_{R_p}-\lambda_{I_p})\sigma_pI_p+(\lambda_{D_p}-\lambda_{I_p})\mu_pI_p.
\end{array}
\end{equation}

The adjoint system is defined as follows:
\begin{equation}\label{adjoint system}
\begin{array}{l}
\dot{\lambda}_{S_p}(t) = (\lambda_{S_p}-\lambda_{A_p})\beta (\sum_{p=1}^{3} A_p+\sum_{p=1}^{3} I_p)+(\lambda_{S_p}-\lambda_{Q_p})(u_p+\gamma_p);\\

\dot{\lambda}_{Q_p}(t) = 0;\\

\dot{\lambda}_{A_p}(t) = (\lambda_{S_p}-\lambda_{A_p})\beta S_p+(\lambda_{A_p}-\lambda_{I_p})\alpha_p k_p+(\lambda_{A_p}-\lambda_{R_p})(1-\alpha_p)\sigma_p;\\

\dot{\lambda}_{I_p}(t) = -f_p'(I_p)+(\lambda_{S_p}-\lambda_{A_p})\beta S_p+(\lambda_{I_p}-\lambda_{R_p})\sigma_p+(\lambda_{I_p}-\lambda_{D_p})\mu_p;\\

\dot{\lambda}_{R_p}(t)=g_p'(R_p);\\

\dot{\lambda}_{D_p}(t) = 0,\\
\end{array}
\end{equation}
with the transversality conditions given by
\begin{equation}\label{transversality}
\lambda_{S_p}(T)=\lambda_{Q_p}(T)=\lambda_{A_p}(T)=\lambda_{I_p}(T)=\lambda_{R_p}(T)=\lambda_{D_p}(T)=0.
\end{equation}

According to Pontryagin's  maximum principle, there exist continuous and piece-wise
continuously differentiable co-state functions $\lambda_r(t),\ r\in\{S_p,Q_p,A_p,I_p,R_p,D_p\}$, $p=1,2,3$  that satisfy (\ref{adjoint system}) and (\ref{transversality}) for  $t\in [0, T]$ together with continuous functions $u^*_1(t)$, $u^*_2(t)$, and $u^*_3(t)$:
\begin{equation}\label{maximum_principle}\begin{array}{l}
(u^*_1, u^*_2, u^*_3)\in \textrm{arg} \min\limits_{{u}_1, {u}_2,{u}_3\in
[0,u_{max}]} H(\lambda_r,S_p,Q_p,A_p,I_p,R_p,D_p,{u}_1,
{u}_2,{u}_3).\end{array}
\end{equation}

Let us define the functions $\varphi_p(t)$ as follows:
\begin{equation}\label{Hamiltonian1}
\begin{array}{l}
\varphi_1(t)=(\lambda_{S_1}(t)-\lambda_{Q_1}(t)), \\
\varphi_2(t)=(\lambda_{S_2}(t)-\lambda_{Q_2}(t)), \\
\varphi_3(t)=(\lambda_{S_3}(t)-\lambda_{Q_3}(t)).
\end{array}
\end{equation}

Let $h_p(\cdot)$ be a concave functions ($h''_p(\cdot)<0$), then according to (\ref{Hamiltonian}), the Hamiltonian is a concave function of $u_p,\ p=\overline{1,3}$. There are two different options for $u_p\in[0,1]$ that minimize the Hamiltonian, i.e., if at the time moment $t$
$$h_p(0)-\varphi_p(t)\cdot 0<h_p(u_{max})-\varphi_p(t)u_{max},$$
or
$$h_p(u_{max})>\varphi_p(t)u_{max},$$
then optimal control is $u_p=0$ (see Fig. 2 (left)); otherwise $u_p=u_{max}$ (see Fig. 2 (right)).

\begin{figure}[ht!]
\begin{center}
  \includegraphics[width=78mm]{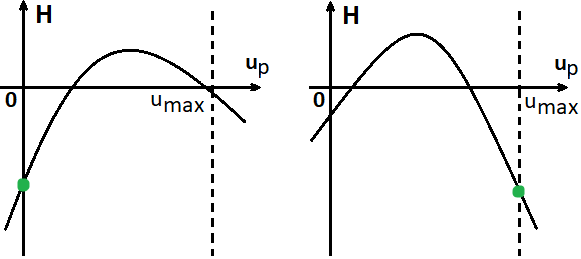}
  \caption{Hamiltonian if functions $h_i(\cdot)$ are concave.}
  \end{center}\label{ham_concave}
\end{figure}

For $p=\overline{1,3}$, the optimal control parameters $u_p(t)$ are defined as follows:
\begin{equation}\label{control}
u^*_p(t)=\left\{\begin{array}{l} 0,\hskip 21pt \mbox{for}\ \ \varphi_p(t)u_{max}<h_p(u_{max}),\\
u_{max}, \ \mbox{for}\ \ \varphi_p(t)u_{max}\geq h_p(u_{max}).\\
\end{array}\right.
\end{equation}

Concave cost functions indeed imply a concave Hamiltonian. The first order conditions for minimization (i.e. the derivatives of the Hamiltonian with respect to the controls equal to zero) are only a necessary - but not sufficient - condition for minimization, being the second requirement a negative second derivative (concavity). If the second derivative is positive, like in this case, it means that the point obtained is a maximum, not a minimum, therefore the solutions offered by this problems are ``corner'' solution, which in this case it means total lockdown or no lockdown at all, as in Figure 2.

\subsection{Functions $h_p(\cdot)$ and the Hamiltonian are strictly convex}
Let $h_p(\cdot)$ be a strictly convex functions ($h''_p(\cdot)>0$), then Hamiltonian is a convex function. Consider the following derivative:
{\begin{equation}\label{der}\begin{array}{l}
\frac{\partial}{\partial x}(h_p(x)-\varphi_p x)\mid_{x=x_p}=0,
\end{array}\end{equation}}
where $x_p\in [0,u_{max}]$, $u^*_p(t)=x_p$. There are three different types of points at which the Hamiltonian reaches its minimum (Fig. 3).
To find them, we need to consider the derivatives of the Hamiltonian at $u_p=0$ and $u_p=u_{max}$. If the derivatives (\ref{der}) at $u_p=0$ are increasing ($h'_p(0)-\varphi_p\ge0$), then the value of the control that minimizes the Hamiltonian is less than 0, and according to our restrictions ($u_p\in [0,u_{max}]$) optimal control will be equal to 0 (Fig. 3a). If the derivatives at $u_p=u_{max}$ are non-increasing ($h'_p(u_{max})-\varphi_p<0$), it means that the value of the control that minimizes the Hamiltonian is greater than $u_{max}$. Hence the optimal control will be $u_{max}$ (Fig. 3c); otherwise, we can find such value $u^*_p\in (0,u_{max})$ (see Fig. 3b).

\begin{figure}[H]
\begin{center}
  \includegraphics[width=105mm]{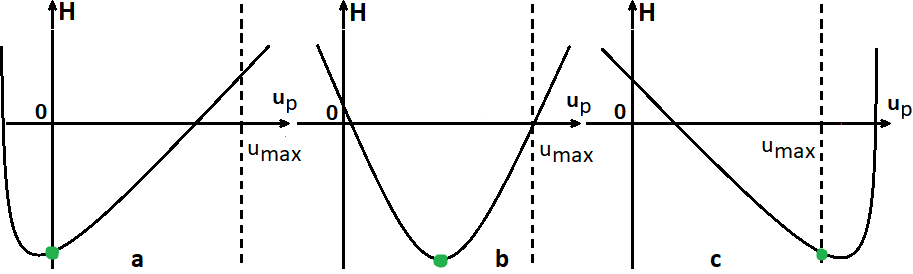}
  \caption{Hamiltonian when functions $h_p(\cdot)$ are convex.}
\end{center}\label{ham_convex}
\end{figure}

Functions $\varphi_p(t)$, $h_p'(t)$, $u^*_p(t)$ are continuous at all $t\in [0, T].$ In this case $h_p$ is strictly convex and $h_p'$ is strictly increasing functions, so $h'_p(0)<h'_p(u_{max})$. Thus there exist  points $t_0$ and $t_1$ $(0<t_0<t_1<T)$ so that conditions (\ref{control1}) are satisfied, and according to $\varphi_p(t)$ are decreasing functions. That is:
\begin{equation}\label{control1}
u^*_p(t)=\left\{\begin{array}{l} 0, \hskip 43pt \mbox{for}\ \varphi_p(t)\leq h'_p(0);\\
h'^{-1}(\varphi_p), \hskip 10pt \mbox{for}\ h'_p(0)<\varphi_p(t)\leq h'_p(u_{max});\\
u_{max}, \hskip 26pt \mbox{for}\ h'_p(u_{max})<\varphi_p(t).
\end{array}\right.
\end{equation}

This means that if the economy presents convex socio-economic costs, that is to say strictly increasing in its production levels and very damaging because of to the lockdown times, then the optimal lockdown time ($t_0$ and $t_1$, where $0<t_0<t_1<T$) presents three cases, i.e., either a nullity  or unit level, but also the case in which the lockdown implies that socio-economic costs take a value between 0 and $u_{max}$.

\bigskip

\newpage

\section{Numerical simulations}

Let us now corroborate our results, performing some numerical simulations.
In this section, we first simulate the economic impact of the pandemic without any policy aimed at limiting this effect. We then simulate the same impact, but assuming the existence of a lockdown policy, which is differentiated between demographic groups. Finally, we simulate the same economic impact of the pandemic, assuming that there is no differentiation between demographic groups, so the same policy is applied to the whole population irrespective to the differentiated effects of the virus among individuals of different ages.
With respect to the theoretical model introduced in the previous section, to simplify and without loss of generality, we simulate only the case of convex cost functions,\footnote{This assumption is based on the fact that in pandemic times, there are not properly economies of scale. For example, Glover et al. (2020) build a model in which economic activity (losses) and disease progression are jointly determined. Palomino et al. (2020) pointed out the high economic costs in the labor market during the pandemic, they develop a Lockdown Working Ability index estimating the potential wage loss and finding that there would be substantial and uneven wage losses across the board; inequality within countries will worsen, as it will between countries. Adams-Prassl et al. (2020) using real-time survey evidence from the UK, US and Germany show that the labor market impacts of COVID-19 differ considerably across countries. Workers in alternative work arrangements and in occupations in which only a small share of tasks can be done from home, are more likely to have reduced their hours, lost their jobs and suffered falls in earnings.  } and moreover, we assume that the profits from recovery are zero since people do not really get any additional benefit to the previous social standard they had before contracting the virus (Acemoglu, 2020).

 In all the simulations, population is normalized to $N=1$ and a whole period of epidemic equal to $T=365$ days. We consider a fixed time period because we assume that at time $T=365$, a vaccine becomes available and then there is no need to implement lockdown policies anymore. We consider as a benchmark example the case of Italy, which consists of 49,581,000 inhabitants. In the remainder of the paper, we will indicate the population with the letter $Z$. The distribution of the population is assumed to be the following: 44.68\% of the people are aged 20-49 (we call this strata as youngs), 27.22\% of the population is aged between 50 and 64 (adults) and the old ones are those aged 65 years or more and represent the 28,10\% of the population.  
As in Acemoglu et al. (2020), we consider a fatality rate by age group, see Table \ref{fatalityrate}

\begin{center}
\begin{table} [hbt!]
    \centering
    \begin{tabular}{c c}
    \hline
        Age group & Fatality rate  \\
        \hline
        20 - 49 & 0.001 \\
        50 - 64 & 0.01 \\
        65+ & 0.06 \\
        \hline
    \end{tabular}
    \caption{Infection Fatality Rate from COVID-19}
    \label{fatalityrate}
\end{table}
\end{center}

Not all the infected develop symptoms: as in Verity et al., (2020) we assume that, on average, only two-third of the cases are sufficiently symptomatic to self-isolate, with a probability for the youngs of developing symptoms equal to $\alpha_1=1/2$, for adults equal to $\alpha_2=2/3$ and for the elderly equal to $\alpha_3=5/6$. 

In the following subsections, we will introduce the methods for calibrating the costs and the benefits functions of our model.

\subsection{Direct costs of infection.}
The direct socioeconomic costs of the infection\footnote{As in Acemoglu et al. (2020), we did not consider for the analysis of costs the strata of population aged less than 20 years.} - as introduced in section \ref{costbenefits} - are represented by the functions: $$Z \cdot f_p(I_p(t)) = Z \cdot E_p^f I_p(t),$$
where $E_p^f$ denote the daily cost of treatment for one infected person belonging to the demographic group $p$. Let $I_p(t)$ be - as usual - the fraction of infected people of group $p$ over the population at all, denoted by $N$. The total cost beared by the Government for the infected belonging to group $p$ is therefore equal to 49,581,000 $\cdot$ $E_p^f I_p(t)$ .
To compute the cost of treatment of the infected $E_p^f$, we must realize that some individuals that are hospitalized need intensive care, but not all. As in Ferguson et al. (2020), we assume a probability of hospitalization of 2.8\% for the youngs, 10\% for the adults and 26\% for the elderly, and those hospitalized who require critical care are 5\% for the youngs, 20\% for the adults and 50\% for the elderly.
According to the most recent medical literature (Gaythorpe et al., 2020), the average stay in hospital is approximately 15 days: 5 days for the young ($\sigma_1=0.2$), 15 days for adults ($\sigma_2=0.066$) and 25 days for the elderly ($\sigma_3=0.04$), with a daily cost per patient of euro 1,500.00 in critical care, and a cost of approximately euro 300.00 in ordinary care.\footnote{For an estimation of the costs of hospitalization in ordinary care for a COVID-19 patient in Italy see  \url{http://www.sossanita.org/archives/9862}, and in critical care see \url{https://www.proiezionidiborsa.it/quanto-ci-costa-un-malato-di-coronavirus/}.}
After discharge, we will assume that the patient is recovered.
To calibrate correctly our model, we therefore estimate, for each age-group, the cost of treatment for the infection by multiplying the probability of hospitalization by its cost. Since not all patients hospitalized need critical care, we ponder the cost of treatment with respect to the probability to need intensive care, which is five times more costly than ordinary care.
Basically,

\bigskip

\begin{center}
    $E_p^f$ = Prob. of hospitalization $\cdot$ (\euro 1500 $\cdot$ prob. of critical care + \euro 300 $\cdot$ (1 - prob. of critical care)).
\end{center}

\bigskip

\noindent which means, for each age group, a daily cost of treatment which is equal to:

\begin{itemize}

\item[-] $E_1^f$= 0.028 $\cdot$(1,500.00 $\cdot$ 0.005 + 300.00 $\cdot$ 0.995) = 8.57,

\item[-] $E_2^f$= 0.1 $\cdot$ (1,500.00 $\cdot$ 0.2 + 300.00 $\cdot$ 0.8) = 54.00,

\item[-] $E_3^f$= 0.26 $\cdot$ (1,500.00 $\cdot$ 0.5 + 300.00 $\cdot$0.5) = 234.00.
\end{itemize}

Obviously, the direct costs of the infection are related only to the Infected group, since the other groups (Susceptibles, Asyntomatic and Quarantined) are assumed not to have symptoms and therefore they do not need to be assisted/cured.

\subsection{Indirect Costs of Infection}

The indirect socioeconomic costs of infection refer basically to the job opportunities lost due to isolation, irrespective of whether it happened as a consequence of an infection or of a lockdown or even a death.
They are group-specific, in the sense that they vary according to the group (infected, quarantined, asymptomatic and dead) and to the age of the individuals, and are related basically to the loss of GDP suffered as a consequence of such events.

Let us begin by illustrating the indirect socioeconomic costs for the group of the infected. An infected person indeed, during the period he must be isolated or treated, of course cannot be able to work. From the Italian Bureau of Statistics (ISTAT) the average yearly nominal GDP per worker in year 2019 was 70,105.30, which leads to a daily income per worker (considering 365 working days in a year) of \euro 192.07.
Following Acemoglu et al. (2020), we assume that the loss of productivity of the young and the adults is the same in this model, while that for the elderly is 26\% of the young (or adult), reflecting the fact that only those aged over 65 that actually work are only 20\% of the total.

These costs can be summarized by the following equation:

\begin{equation}
    Z \cdot h^h_p(I_p(t)) = Z \cdot E^I_pI_p(t).
\end{equation}

where $E^I_p$ is equal to:
\begin{itemize}
    \item[-] $E_1^I$=192.07,
    \item[-] $E_2^I$= 192.07,
    \item[-] $E_3^I= 192.07\cdot 0.26= 49.94$.
\end{itemize}

The indirect socioeconomic costs for the individuals quarantined by by law (that is to say, the fraction $u_p$ of the susceptible) can be defined as:

\begin{equation}
    Z \cdot h^h_p(u_p(t), \gamma_p) = Z \cdot E^S_p(u_p(t)+\gamma_p)^2S_p(t).
\end{equation}


We remind the fact that the government treats all the infected she/he is able to detect. Treatment may consists in hospitalisation if the symptoms are severe, or isolation if the individual suffers mild symptoms. The policy of lockdown is therefore directed to the susceptible only, which, not being sick, may decide to act freely as they did before the pandemics. We assume that these costs grow exponentially with respect to the number of infected individuals, because we assume diminishing return to scale for the labor production factor (as the number of infected individual increases, the productivity of the non-isolated individual decreases, too, up to a point beyond which additional increases in output cannot be made).

Each individual that is isolated experiences a productivity loss. We assumed previously that an infected individual, for each day he stays isolated or hospitalized, experiences a productivity loss of 100\% of his daily income. A susceptible that is auto-isolated or quarantined by law instead, may keep his ordinary activity (whenever possible) from home, so his productivity loss is lower.

Following Acemoglu et al. (2020), the susceptibles that are auto-isolated or quarantined have a loss in productivity of 70\%, that is to say that $E_p^S$, are represented by the following values:

\begin{itemize}
    \item[-] $E_1^S$=192.07 $\cdot$ 0.7 = 134.45,
    \item[-] $E_2^S$= 192.07 $\cdot$ 0.7 = 134.45,
    \item[-] $E_3^S$= 49.94 $\cdot$ 0.7 = 34.96.
\end{itemize}

Finally, each person that, as a consequence of the infection, dies, produces a loss for the society which in principle is very difficult to estimate.
The function for the indirect costs of a death is the following:
\begin{equation}
    Z \cdot h^h_p(D_p) = Z \cdot E^D_pD_p(t).
\end{equation}

In our paper, we will assume that the cost of a life is, for the three groups:

\begin{itemize}
    \item[-] $E_1^D$= 2,800,000.00,
    \item[-] $E_2^D$= 2,000,000.00,
    \item[-] $E_3^D$= 273,000.00.
\end{itemize}

\bigskip

\subsection{Numerical results}

In the previous two subsections we illustrated how we calibrated the model for numerical simulations.\footnote{For everything is not explained or mentioned, the reader can refer to tables A.1 and A.2 in the appendix.}
As already pointed out, our numerical simulations refer to the case where the cost functions $h_p$ are strictly convex, and there are no profits from recovery.

We will perform two experiments: the first, where a policy determined by optimal control is compared with no policy. An estimation of the dynamics of the groups of Susceptibles, Quarantined, Asymptomatic, Infected, Recovered and Dead is performed through time. Moreover, we computed the costs associated to these two policies (an optimal policy as compared to \emph{laissez faire}).
The second experiment simulate an optimal uniform policy, without distinguishing between demographic groups. This experiment is important to estimate the costs associated to the optimal policy and compare it with the optimal targeted policy performed in experiment 1.


\bigskip

\textbf{ Experiment 1.} Using data of column 1 in table A.1 (see the appendix), we estimated the behaviour through time of the groups of Susceptible, Quarantined, Asymptomatic, Infected, Recovered and Dead, for the three demographic groups. Figure 4 and 5 describe the behaviour of these groups without any policy aimed at limiting the spread of the virus, while figures 6 - 9 describe the same variables when this policy is introduced. For this experiment and the next one, we considered a population of 49,581,000 people, with the following demographic strata: 22,150,000 aged 20-49 (44.68\%); 13,498,000 aged 50-64 (27.22\%), and 13,933,000 aged 65+ (28.10\%). As it is possible to notice, in the absence of a policy, the total number of infected increases up to the the peak of day 35 (about) and decreases thereafter. The total number of infected (and asymptomatic) during the whole period of pandemic may reach about 21 millions (42,4\% of the whole population), if no policy is implemented. Conversely, when a targeted policy (based on demographic groups) is undertaken, which is aimed at reducing the contagion, a different scenario appears. The total number of infected and asymptomatic indeed more than halves; from about to 21m to about 10m at the peak of the pandemic, while almost half of the population is quarantined (see figure 6). Moreover, the peak of the pandemic occurs after 20-22 days from the beginning, against 35 days in an uncontrolled case.

\begin{figure}[H]
\begin{center}
\includegraphics[width=53mm]{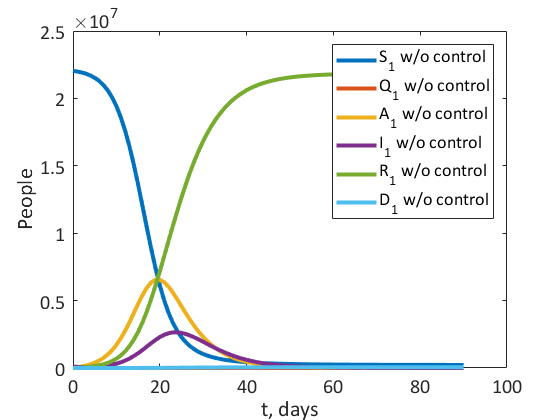}
\includegraphics[width=53mm]{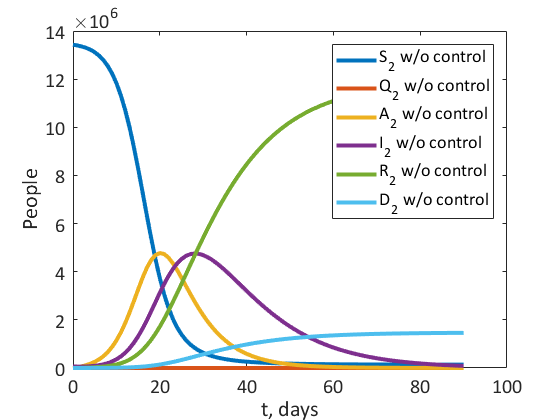}
\includegraphics[width=53mm]{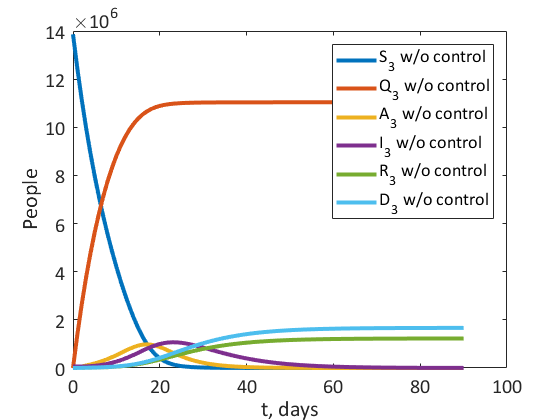}
\caption{\footnotesize Experiment 1. The behaviour of the system in Young (left), Adult (center), and Old (right) subpopulations in uncontrolled case.}
\end{center}\label{}
\end{figure}

\begin{figure}[H]
\begin{center}
\includegraphics[width=53mm]{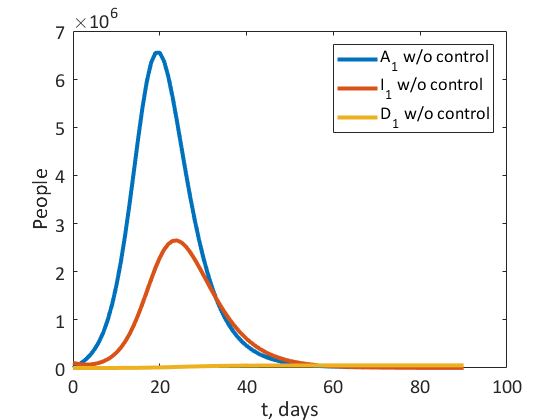}
\includegraphics[width=53mm]{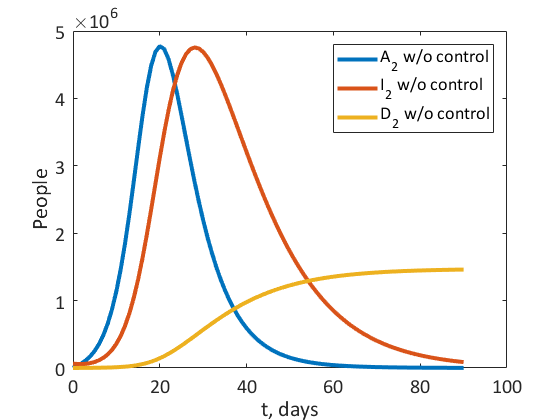}
\includegraphics[width=53mm]{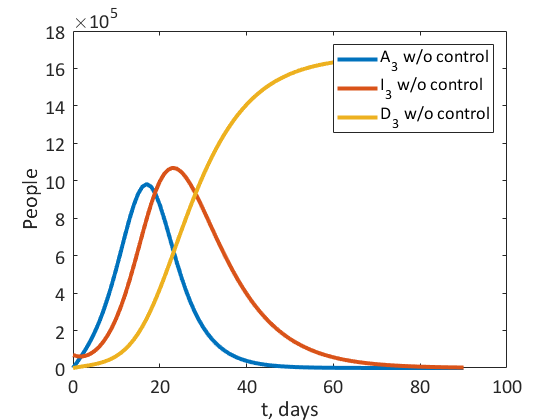}
\caption{\footnotesize Experiment 1. Number of Asymptomatic, Infected, and Dead in uncontrolled case.}
\end{center}\label{}
\end{figure}

It is obviously clear that the most penalised group in terms of number deaths is the group of the elderly. Without any control, indeed, this group loses more than 1.6 millions of people (about 11,5\%), against 1.5 millions losses in the adult group (11,1\%) and a negligible amount among the youngs. Without control, those who are quarantined are only a fraction of adult and old groups that arbitrarily decide to self-quarantine (0.3\% of the adults, and 10\% of the elderly, as it is possible to see in the table in the appendix, column ``Exp. 1''.)

\begin{figure}[H]

\begin{center}
\includegraphics[width=53mm]{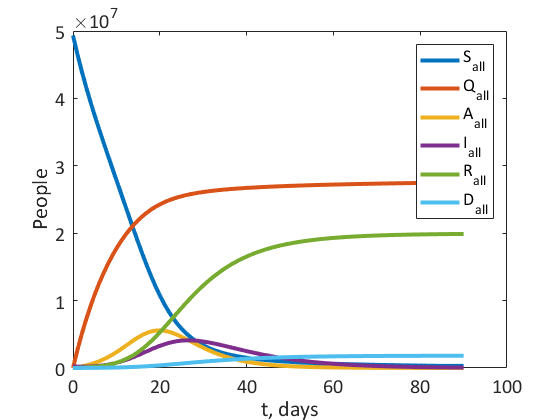}
\caption{\footnotesize Experiment 1. Sum of fractions of all three groups (young, adult, and old) in controlled case.}
\end{center}\label{Compare}
\end{figure}

When a control policy is implemented, globally about a half population is quarantined. Figure 6 display a total of 25 millions people isolated (which, considering that the total population of this experiment amounts to 49.5 millions, it's about one half).

\begin{figure}[H]
\begin{center}
\includegraphics[width=53mm]{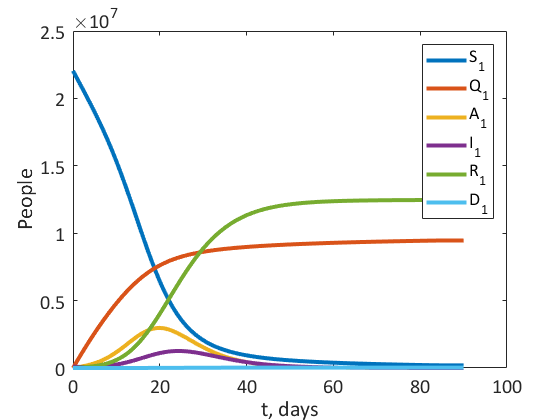}
\includegraphics[width=53mm]{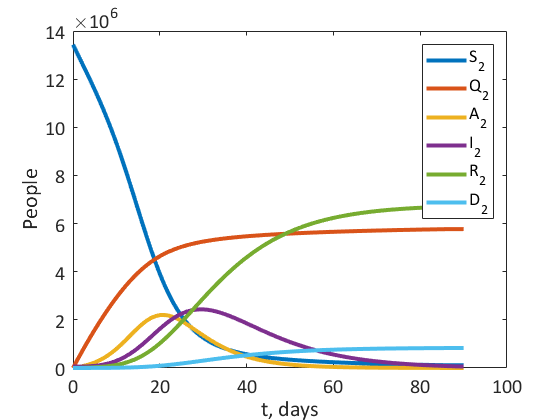}
\includegraphics[width=53mm]{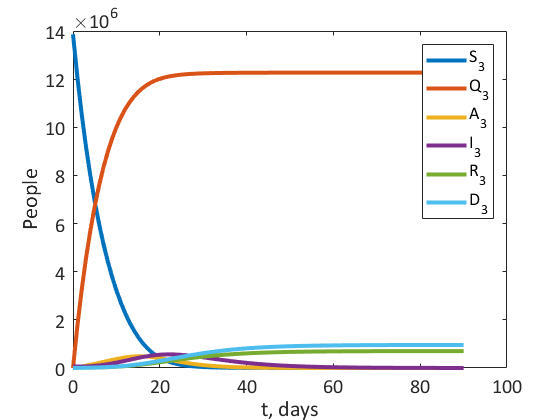}
\caption{\footnotesize Experiment 1. The behaviour of the system in Young (left), Adult (center), and Old (right) subpopulations  in controlled case.}
\end{center}\label{}
\end{figure}

Figure 7 displays the behaviour of the groups of the Susceptibles, Quarantined, Asymptomatic, Infected, Recovered and Dead for the three demographic groups of youngs, adults and old.
Intuitively, our simulations show that about 10 millions of young people (which corresponds to about 45\%) are quarantined at the end of the pandemic, against 6 millions of adults (44.4\%) and 12 millions of old people (86\%). Quarantined people include those who are quarantined by law and those who deliberately quarantined themselves. The reason for this effect relies on the fact that a young person has a higher probability to be asymptomatic and affects other people belonging to the other two groups. Figure 8 highlights better this phenomena.

\begin{figure}[H]
\begin{center}
\includegraphics[width=53mm]{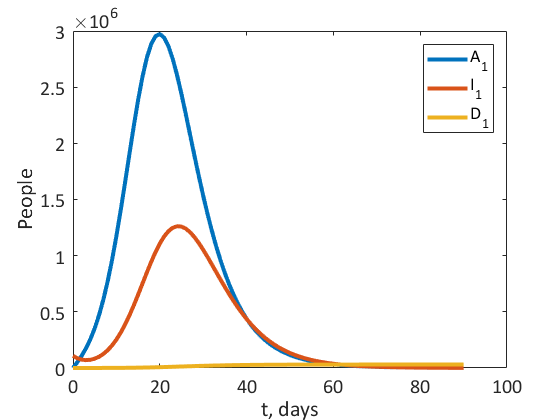}
\includegraphics[width=53mm]{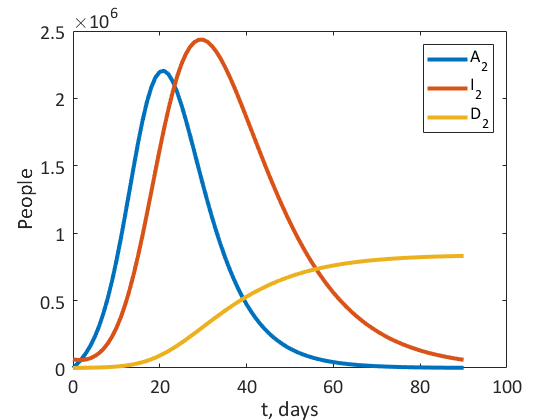}
\includegraphics[width=53mm]{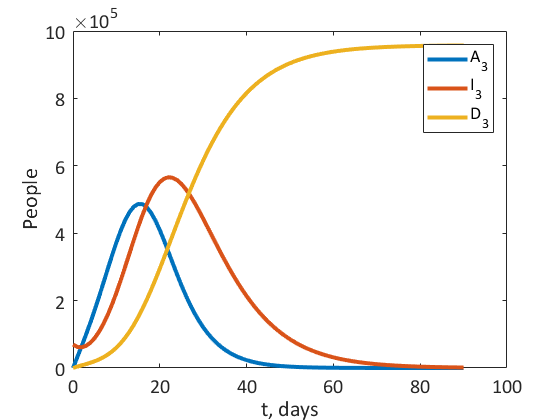}
\caption{\footnotesize Experiment 1. Number of Asymptomatic, Infected, and Dead in controlled case.}
\end{center}\label{}
\end{figure}

Since the transmission rate of the virus is the same among the three groups, it is clear that the group who is at higher risk of not developing symptoms (so not being detected) has to be quarantined in order to avoid deaths in the other two groups, despite the cost of a life lost is unbalanced/skewed towards the youngs and the adults.

An optimal control policy here assumes at the beginning of the pandemics (at t=0), an immediate isolation of 5.5 millions of youngs (24.83\%), and 3.3 millions of adults and old people (which correspond to 24.4 and 23.7\%, respectively). The subsequent number of people isolated decrease day by day, by the way (figure 9, on the left, shows each day of the pandemic the number of isolated person, for each demographic group. Once a person is isolated, it is assumed to be quarantined and stays in isolation up to the end of the pandemic).

\begin{figure}[H]
\begin{center}
\includegraphics[width=53mm]{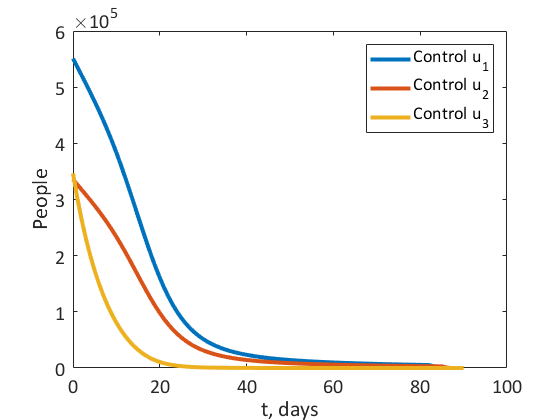}
\includegraphics[width=53mm]{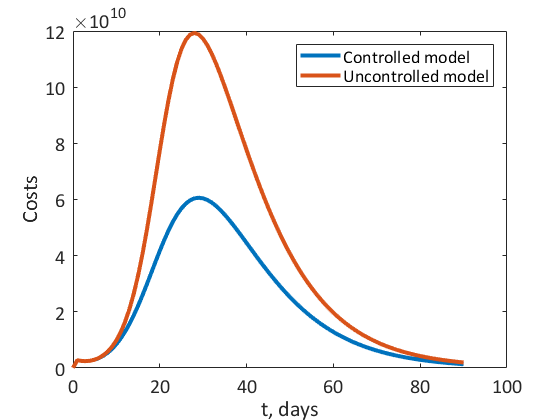}
\caption{\footnotesize Experiment 1. Structure of the optimal control (left). Aggregated system costs are $J=2.0418\cdot 10^{12}$ m.u. in controlled case and $J=3.5809\cdot 10^{12}$ m.u. in uncontrolled case (right).}
\end{center}\label{}
\end{figure}

The aggregated costs of the no-intervention policy amount to \euro 3.5809$\cdot 10^{12}$, against a control policy which generates costs for \euro 2.018$\cdot 10^{12}$ (see figure 9, on the right.)

\bigskip

\textbf{ Experiment 2.}
In this second experiment, we run the simulation using data from the first experiment (and averaging when necessary) and considering the whole population as a single group. This experiment was performed, basically, to compare the aggregated costs arising when the policy implemented to limit the spread of the virus does not exploit the demographic differences (and therefore the differences in the health outcomes generated by the virus) between groups.


\begin{figure}[H]
\begin{center}
\includegraphics[width=53mm]{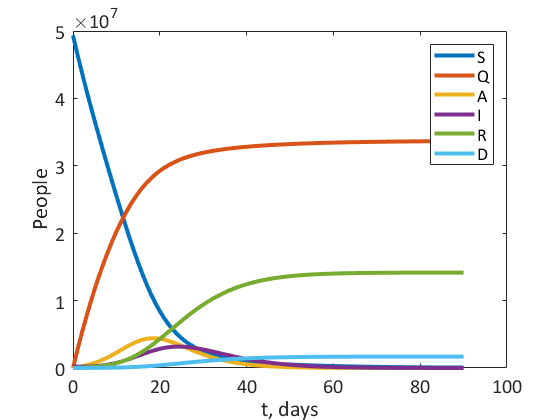}
\includegraphics[width=53mm]{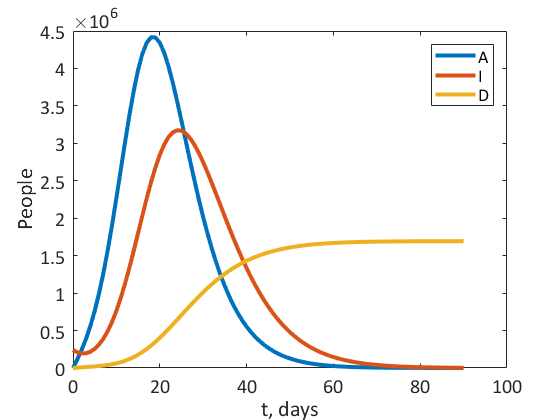}
\caption{\footnotesize Experiment 2. The behaviour of the system when the whole population is considered as a one group in controlled case. (left) Number of Asymptomatic, Infected, and Dead in controlled case. (right)}
\end{center}\label{Compare}
\end{figure}

What catches immediately the eyes is that at the end of the pandemics, 33 millions of people are quarantined, that is to say, 66.5\%, against a 50\% of an optimal targeted policy. This of course has its consequences on the aggregated costs: this untargeted policy indeed generates cost for $3.1905\cdot 10^{12}$ against a cost of  $2.0418\cdot 10^{12}$ of a targeted policy.

\begin{figure}[H]
\begin{center}
\includegraphics[width=53mm]{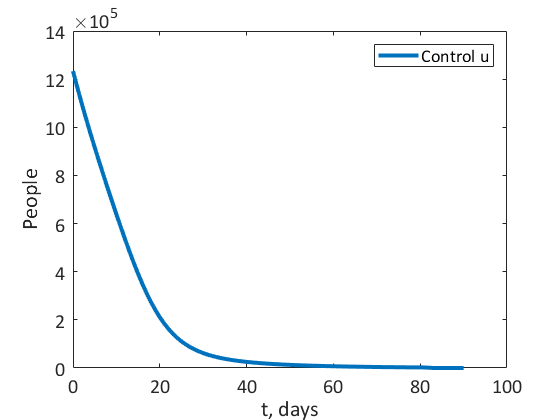}
\includegraphics[width=53mm]{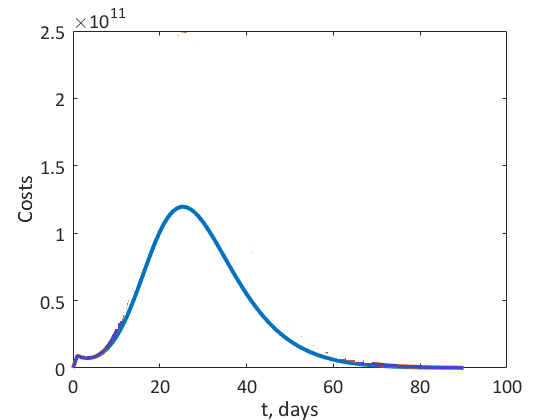}
\caption{\footnotesize Experiment 2. Structure of the optimal control (left). Aggregated system costs are $J=3.1905\cdot 10^{12}$ m.u. in controlled case (right).}
\end{center}\label{}
\end{figure}


It is therefore pretty evident that a selective, targeted policy outperforms a uniform policy in terms of the costs that it generates, because a uniform policy generates 56\% more costs than a targeted one.
We do not observe a substantial different durations of the policies, which last, in both cases, about 50 days.

It is important, then, to calibrate the controls selectively on the different demographic groups so as to exploit the differential impact that the pandemic has on those populations.
Even though this does not allow to shorten the lockdown, but it allows to keep in the job market the more productive and young people (therefore, those who suffer less the consequences of the virus) thus limiting the costs associated to the policy of lockdown.

\newpage
\section{Concluding remarks}

In this paper, we modified a classical SIR model in order to develop an optimal control policy of lockdown aimed at minimizing a cost-profit social function.

Our main results are that different lockdown policies (total lockdown, \emph{laissez faire} and a mixed policy, that is to say, partial lockdown) may be explained as a consequence of different cost structures of a country.

If, indeed, the two corner solutions, that is to say, total lokdown and no lockdown may arise in the presence of concave cost functions, under the hypothesis of convex cost functions an internal solution is optimal.

Under the hypothesis of convex cost function, that is to say, when it is optimal to apply mixed strategy (i.e. a partial lockdown), we also showed that a targeted lockdown policy based on demographic groups outperforms a uniform policy based on the whole population, because, in according to our simulations, it allows to save about 56\% of the costs.

From a political perspective, it is important for policymakers that have to manage an emergency like those of the recent COVID-19 pandemic in the best way.
This means to implement an effective policy of lockdown (in the absence of any medical treatment like a vaccine that can help develop immunity for the whole population) which produces the minor costs for the society.

A lockdown indeed has a deep impact on the families' economies (especially of those whose income is linked to activities that are at greater risk of contagion), but catching the virus may be tragic as well, because the long-run consequences may be serious and responsible for income losses. A responsible management of such health emergency is of crucial importance for the well being of a society, and limits subsequent hard economic policies aimed at restoring the basic subsistence levels of the population.

\bigskip

Finally, we must point out that our optimal policy does not consider additional prescriptions like social distancing, facial mask wearing and so on. Our results indeed may be widened by imposing, for those who are not involved by any lockdown, of such additionally policies, which are assumed to reduce the propagation rate of the virus, according to the medical literature.

Moreover, we consider only age as discriminant for a given policy. COVID-19 has proven to have differentiated impacts on the health of the individuals depending on previous comorbidities, given the age. Our results could be improved - and better policies implemented - if, instead of considering just three groups based on age, we could distinguish groups based on pregressed comorbidities and age. We decided however to leave these experiments to future research.

Not less important is to let the reader notice that despite our effort to calibrate the model in the most realistic way, there is a huge uncertainty about the parameters that govern the dynamics of the infection, and therefore there is a potential bias in the determination of the costs associated to the policies undertaken, irrespective of whether this policy implies a lockdown or not.

\bigskip

We wants to make our apologies to those who, reading this paper, have felt indignant by the aseptic presentation of the argument, which has been focused in a mere computation of economic costs and profits of a policy intervention, without regards to the psychological impact that lives losts may have for the relatives and the society at all.
We are aware that the pandemic has caused a loss many human lifes, and a dissertation of the argument that is focused exclusively on costs and profits could hurt the feelings of many people.
As economists, we are trained since the beginning of our studies to think in those terms, so we have a sort of ``mental deviation'' that induces us to analyze strategic situations like this one in an impersonal, aseptic and distant way, on the basis of economic indicators only.
In spite of everything, we believe that this kind of presentation could be useful for policymakers who wants to manage such health emergency in an effective (and possibly efficient) way.

\newpage

\section*{Appendix}


\hskip -20pt \begin{center}
\footnotesize
\begin{tabular}{| l | c | c |   }
\hline
\multicolumn{3}{|l |}{ \textbf{Table A.1.: Parameters used for simulations in the experiments} } \\
\hline\hline
Parameters' name & Exp. 1 & Exp. 2   \\
\hline

Population  & 49,581,000 & 49,581,000    \\
\hline
Fraction of & 0.4446 & 0.9951    \\
Susceptibles & 0.2709 &     \\
at time 0 (\textbf{$S_p^0$})    & 0.2796 &     \\
\hline
Fraction of & 0 & 0    \\
Quarantined & 0 &    \\
at time 0 (\textbf{$Q_p^0$})    & 0 &     \\
\hline
Fraction of & 0 & 0    \\
Asymptomatic & 0 &    \\
at time 0 (\textbf{$A_p^0$})    & 0 &     \\
\hline
Fraction of & 0.0022 & 0.0049     \\
Infected     & 0.0013 &       \\
at time 0 (\textbf{$I_p^0$})    & 0.0014 &     \\
\hline
Fraction of & 0 & 0     \\
Recovered & 0 & 0    \\
at time 0 (\textbf{$R_p^0$})    & 0 & 0    \\
\hline
Fraction of & 0 & 0    \\
Dead     & 0 & 0   \\
at time 0 (\textbf{$D_p^0$})    & 0 & 0  \\
\hline
Recovery rate (\textbf{$\sigma_p$}) & 0.2 & 0.118   \\
                   & 0.066 &        \\
                   & 0.04 &      \\
                   \hline

Transition rate & 0 & 0.028  \\
from $S_p$ to $Q_p$ \textbf{$(\gamma_p)$} & 0 &    \\
                & 0.1 &    \\
                \hline
Death rates (\textbf{$\mu_p$}) & 0.001 & 0.02  \\
            & 0.01 &    \\
            & 0.06 &    \\
            \hline
Infection rate & 0.48  & 0.48    \\
from $S_p$ to $I_p$  (\textbf{$\beta$})&     &      \\
\hline

Asymptomatic to & 0.1923 & 0.1923     \\
Infected (\textbf{$k_p$}) & & \\
\hline
Probability that  & 0.5 & 0.636       \\
node will have             & 0.66 &       \\
symptoms (\textbf{$\alpha_p$})            & 0.83 &       \\
            \hline

Costs of treatment \textbf{$E^f_p$} & 200.64 & 236.413  \\
            & 246.07 &      \\
            & 283.94 &      \\
            \hline

Costs of quarantine \textbf{$E^h_p$} & 134.45 & 106.493   \\
            & 134.45 &    \\
            & 34.96 &    \\

            \hline
Cost of death \textbf{$E^D_p$} & 2800000 & 1872153   \\
            & 2000000 &    \\
            & 273000 &    \\
            \hline
Aggregated costs \textbf{$J$} &  $3.5809\cdot 10^{12}$ & $5.9475\cdot 10^{12}$  \\
 Uncontrolled case & & \\
\hline
Aggregated costs \textbf{$J$} &  $2.0418\cdot 10^{12}$ & $3.1905\cdot 10^{12}$   \\
 Controlled case & &\\
\hline
\end{tabular}
\end{center}

\begin{center}
\begin{tabular}{| l | c | c | }
\hline
\multicolumn{3}{|l |}{ \textbf{Table A.2.: Maximum values} } \\
\hline\hline
Max value of $I_p$ &  1262100 & 3175100  \\
 &  2435100 &   \\
 &  566020 &   \\
\hline
Max value of $A_p$ & 2972500 & 4416100      \\
 &  2204600 &      \\
  &  486640 &      \\
\hline
Max value of $D_p$ & 30750 & 1693400     \\
 &  831560 &     \\
  &  957180 &   \\
\hline

Max value of $I_p$  & 2648700 & 6574100   \\
in uncontrolled case & 4755800 &      \\
 & 1069000  &     \\
\hline

Max value of $A_p$ & 6549700 &  9644000 \\
in uncontrolled case & 4775000 &    \\
 & 982260  &     \\
\hline

Max value of $D_p$ & 53738 & 3156800     \\
in uncontrolled case & 1460500 &     \\
 & 1665600  &    \\
\hline

\end{tabular}
\end{center}

\end{document}